\documentclass{amsart}

\usepackage{amscd,amssymb,graphics}
\usepackage{hyperref}
\usepackage{amsfonts}
\usepackage{amsmath}
\usepackage{enumerate}
\usepackage{tikz}

\input xy
\xyoption{all}

\usepackage{vmargin}
\setmarginsrb{20mm}{10mm}{20mm}{13mm}%
{10mm}{7mm}{10mm}{10mm}

\newtheorem{thm}{Theorem}[section]

\newtheorem{corol}[thm]{Corollary}
\newtheorem{lemma}[thm]{Lemma}
\newtheorem{prop}[thm]{Proposition}

\newtheorem{definition}[thm]{Definition}
\numberwithin{equation}{section}
\theoremstyle{remark}
\newtheorem{remark}[thm]{Remark}
\newtheorem{example}[thm]{Example}

\newcommand{\ben}{\begin{enumerate}}
	\newcommand{\een}{\end{enumerate}}

\def\R {{\mathbb R}}

\def\N{{\mathbb N}}

\def\End{\operatorname{End}}
\def\Aut{\operatorname{Aut}}

\def\AT{\operatorname{AT}}

\newcommand{\UT}{\mathrm{UT}}
\def\Np{{\mathbb N_+}}
\newcommand{\Pfin}[1]{\mathcal P_{\mathrm{fin}}(#1)}

\begin{document}	
	
	\title[The Addition Theorem for the Algebraic Entropy of Nilpotent Groups]{The Addition Theorem for the Algebraic Entropy of Nilpotent Groups}
	\author[]{Menachem Shlossberg}
	\address[M. Shlossberg]
	{\hfill\break School of Computer Science 
		\hfill\break Reichman University 
		\hfill\break 8 Ha'universita Street, Herzliya,  4610101
		\hfill\break Israel}
	\email{menachem.shlossberg@post.runi.ac.il}

\begin{abstract}
	The Addition Theorem for the algebraic entropy of group endomorphisms was first
	proved for torsion abelian groups by Dikranjan, Goldsmith, Salce and
	Zanardo~\cite{DGSZ}, and was later extended to all abelian groups by Dikranjan
	and Giordano Bruno~\cite{DGBabelian}. Shlossberg~\cite{S} subsequently proved
	it for torsion nilpotent groups of class~$2$.
	As our main result, we prove the Addition Theorem for endomorphisms of
	nilpotent groups of arbitrary nilpotency class, without any torsion assumption.
	In the torsion case, this yields in particular that, if $G$ is a torsion
	nilpotent group and $\phi\in\End(G)$, then either $h(\phi)=\infty$ or
	$h(\phi)=\log(\alpha)$ for some positive integer $\alpha$.
	We further prove that, for every group $G$ and every automorphism
	$\phi\in\Aut(G)$, the Addition Theorem holds with respect to each term of the
	upper central series. As a consequence, the Addition Theorem holds for
	automorphisms of every $\omega$-hypercentral group.
	Finally, we establish a reduction principle: if $\mathfrak X$ is a class of
	locally finite groups closed under taking subgroups and quotients, then the
	Addition Theorem holds for endomorphisms of groups in $\mathfrak X$ if and only
	if it holds for endomorphisms of those groups in $\mathfrak X$ generated by
	bounded subsets, that is, subsets whose element orders are uniformly bounded.
\end{abstract}

	\maketitle
\section{Introduction}
Algebraic entropy is a numerical invariant measuring the asymptotic growth
of trajectories of finite subsets under iterates of group endomorphisms.
Since its introduction by Adler, Konheim and McAndrew~\cite{AKM}, algebraic
entropy has been widely studied, especially in the abelian setting, where a
rich and well-understood theory has been developed; see, for example,
\cite{Pet,W,DGSZ,DGBabelian}. In particular, Weiss~\cite{W} and Peters~\cite{Pet}
established Bridge Theorems connecting algebraic entropy with topological
entropy.

A fundamental structural property of algebraic entropy is the
\emph{Addition Theorem} (AT), which describes how entropy behaves with
respect to short exact sequences of groups that are invariant under a
given endomorphism. Roughly speaking, the Addition Theorem expresses the
entropy of an endomorphism as the sum of the entropies induced on a normal
subgroup invariant under that endomorphism and on the corresponding quotient.

In the abelian setting, the Addition Theorem was first proved for torsion
abelian groups by Dikranjan, Goldsmith, Salce and Zanardo~\cite[Theorem~3.1]{DGSZ}, and later
extended to all abelian groups by Dikranjan and Giordano Bruno~\cite{DGBabelian}
(see also~\cite{AKM,Pet,W} for classical origins and early developments).
Nevertheless, Giordano Bruno and Spiga showed in~\cite[Example~2.7]{GBSp} that the Addition
Theorem fails in general for automorphisms of metabelian groups. This
counterexample is closely related to growth phenomena. Indeed, algebraic
entropy measures the exponential growth rate of trajectories of finite subsets
under endomorphisms. For the identity endomorphism, if $F$ is a finite symmetric
generating set containing the identity, then the $n$-th trajectory is $F^n$, the
ball of radius $n$ with respect to $F$. Thus the corresponding entropy records
the exponential growth rate of the group. In this context, nilpotent groups
occupy a distinguished position: finitely generated nilpotent groups have
polynomial growth by a classical theorem of Wolf~\cite{Wolf}, while Gromov's theorem characterizes finitely
generated groups of polynomial growth as virtually nilpotent~\cite{Gromov}.
These results make nilpotent groups a natural class in which to investigate
the validity of the Addition Theorem beyond the abelian setting.

The Addition Theorem was previously known for certain classes of locally finite
non-abelian groups. In particular,
Shlossberg~\cite{S} proved it for endomorphisms of torsion nilpotent groups
of class~$2$; see also \cite{DGSZ,GBS,GBST,XST} for related results in the
locally finite setting.

The purpose of this paper is to remove both the restriction on the
nilpotency class and the torsion assumption. Our main result establishes
the Addition Theorem for endomorphisms of arbitrary nilpotent groups.
Thus, within the class of nilpotent groups, the behavior of algebraic
entropy with respect to invariant normal subgroups and the corresponding
quotients is completely determined by the Addition Theorem.

As a first application, we show that if $G$ is a torsion nilpotent group and
$\phi\in\End(G)$, then either $h(\phi)=\infty$ or
$h(\phi)=\log(\alpha)$ for some $\alpha\in\Np$.
We next prove that, for every group $G$ and every automorphism
$\phi\in\Aut(G)$, the Addition Theorem holds with respect to each term of the
upper central series. Using a continuity argument along
directed unions, we deduce that the Addition Theorem holds for automorphisms
of every $\omega$-hypercentral group. Moreover, when such a group is torsion,
either the entropy of an automorphism is infinite or it is of the form
$\log(\alpha)$ for some $\alpha\in\Np$.

Finally, we establish a reduction principle for classes of locally finite
groups closed under taking subgroups and quotients. For every such class
$\mathfrak X$, the Addition Theorem holds for endomorphisms of all groups
in $\mathfrak X$ if and only if it holds for endomorphisms of those groups in
$\mathfrak X$ that are generated by bounded subsets, that is, subsets whose
element orders are uniformly bounded.

	\section{Preliminaries}\label{sec:preliminaries}
	Following \cite{DG-islam}, we recall the definition of algebraic entropy for
	endomorphisms of arbitrary groups. For a set $X$, let
	\[
	\Pfin{X}
	=
	\{F\subseteq X:0<|F|<\infty\}
	\]
	denote the family of all non-empty finite subsets of $X$.
	
	\begin{remark}\label{rem:fekete}
		We will repeatedly use the following well-known fact: if
		$(a_n)_{n\in\Np}$ is a subadditive sequence of real numbers, that is,
		\[
		a_{n+m}\le a_n+a_m \quad \text{for all } n,m\in\Np,
		\]
		then the limit
		\[
		\lim_{n\to\infty}\frac{a_n}{n}
		\]
		exists and equals $\inf_{n\in\Np}a_n/n$.
		This is the classical lemma of Fekete.
	\end{remark}
	
	\begin{definition}\label{def:entropy}
		Let $G$ be a group, let $\phi\in\End(G)$, and let
		$F\in\Pfin{G}$. For $n\in\Np$, the \emph{$n$-th trajectory of
		$F$ under $\phi$} is
		\[
		T_n(\phi,F):=F\,\phi(F)\cdots\phi^{n-1}(F).
		\]
		The \emph{algebraic entropy of $\phi$ with respect to $F$} is
		\[
		H(\phi,F)
		:=\limsup_{n\to\infty}\frac{1}{n}\log\bigl|T_n(\phi,F)\bigr|.
		\]
		The sequence $\bigl(\log|T_n(\phi,F)|\bigr)_{n\in\Np}$ is
		subadditive, since
		\[
		T_{n+m}(\phi,F)
		=
		T_n(\phi,F)\,\phi^n\bigl(T_m(\phi,F)\bigr)
		\]
		for all $n,m\in\Np$. Hence, by Remark~\ref{rem:fekete}, the above
		$\limsup$ is in fact a limit. The \emph{algebraic entropy} of $\phi$ is
		\[
		h(\phi)
		:=
		\sup_{F\in\Pfin{G}}H(\phi,F).
		\]
	\end{definition}
	
	\begin{remark}\label{rem:locally-finite-cofinal}
		The function $H(\phi,\mathord{-})$ is monotone on $\Pfin{G}$ with
		respect to inclusion. If $G$ is locally finite, every $F\in\Pfin{G}$
		is contained in a finite subgroup of $G$. Hence the finite subgroups form
		a cofinal family in $\Pfin{G}$, and therefore
		\[
		h(\phi)
		=
		\sup\{H(\phi,E):E\leq G\text{ is a finite subgroup}\}.
		\]
	\end{remark}
	\subsection{Notation and terminology}
		The sets of non-negative reals, non-negative integers and positive natural numbers
	are denoted by $\R_{\geq 0}$, $\N$ and $\Np$, respectively.
	Throughout the paper, the notation $H \unlhd G$ indicates that $H$ is a normal subgroup of $G$.	An element $x$ of a group $G$ is \emph{torsion} if the subgroup $\langle x\rangle$ generated by $x$ is finite.
	Recall that a group $G$ is \emph{locally finite} if every finitely generated subgroup of $G$ is finite.
	The \emph{derived series} of a group $G$ is defined by
	\[
	G^{(0)}=G,\qquad G^{(i+1)}=[G^{(i)},G^{(i)}]
	\quad\text{for }i\geq 0.
	\]
	A group $G$ is called \emph{solvable} if $G^{(k)}=\{1\}$ for some
	$k\in\N$. The least such $k$ is called the \emph{derived length} of
	$G$. In particular, $G$ is \emph{metabelian} if
	$G''=G^{(2)}=\{1\}$.
	We also write $\gamma_1(G)=G$, $\gamma_{i+1}(G)=[\gamma_i(G),G]$ for the lower
	central series. In particular, $G'=\gamma_2(G)$ and if $G$ is \emph{nilpotent}
	of class at most $n$, then $\gamma_{n+1}(G)=\{1\}$ and $\gamma_n(G)\le Z(G)$,
	where $Z(G)$ denotes the \emph{center} of $G$.
	
	A group $G$ is said to be \emph{locally solvable} (respectively, \emph{locally
		nilpotent}) if every finitely generated subgroup of $G$ is solvable
	(respectively, nilpotent).
	It is well-known that every torsion locally solvable group is locally finite
	(see \cite[Proposition 1.1.5]{DX}).
	
	We denote by $\End(G)$ the set of all endomorphisms of $G$,
	and by $\Aut(G)$ its subset of automorphisms.
	
	If $\phi\in \End(G)$, then a subgroup $H$ of $G$ is called \emph{$\phi$-invariant}
	if $\phi(H)\subseteq H$, and $H$ is \emph{$\phi$-stable} if $\phi\in \Aut(G)$ and $\phi(H)=H$.
	A subgroup $H$ of $G$ is \emph{characteristic} if $H$ is $\phi$-invariant for every
	$\phi\in \Aut(G)$, and \emph{fully characteristic} if the same holds for every
	$\phi\in \End(G)$. For every $n\in\Np$ we let $G[n]$ be the fully characteristic subgroup
	of $G$ generated by the set $\{x\in G : x^n=1\}$.
	
	We recall that the \emph{upper central series} of a group $G$ is defined inductively
	by $Z_0(G)=\{1\}$ and
	\[
	Z_{n+1}(G)/Z_n(G)=Z\bigl(G/Z_n(G)\bigr)\qquad(n\in\N).
	\]
	Each $Z_n(G)$ is a characteristic subgroup of $G$.
	A group $G$ is said to be \emph{$\omega$-hypercentral} if
	\[
	G=\bigcup_{n\in\N} Z_n(G).
	\]
	Recall that hypercentral groups (in the transfinite sense) are locally
	nilpotent by a classical result of Mal'cev (see, e.g., \cite[p.~8]{DX}).
	In particular, every $\omega$-hypercentral group is locally nilpotent.

	\section{Auxiliary results}\label{sec:reductions}
	We fix the following terminology for the Addition Theorem in a class of
	groups, which will be used throughout the paper:
	\begin{definition}
		Let $\mathfrak X$ be a class of  groups closed under taking subgroups and quotients.\ben \item	We say that $\AT(G,\phi,H)$ holds for a group $G\in \mathfrak X$, $\phi\in \End(G)$ and a $\phi$-invariant normal subgroup $H$ of $G$ if
		\[h(\phi) = h(\phi \upharpoonright_H)+h (\bar \phi),\] where $\bar \phi=\bar\phi_{G/H}\in \End(G/H)$ is the map induced by $\phi.$ \item  The  Addition Theorem  holds in $\mathfrak X$  for endomorphisms if $\AT(G,\phi,H)$ holds for every $G\in \mathfrak X, \ \phi\in \End(G)$  and every $\phi$-invariant normal subgroup $H$ of $G.$
		\item The  Addition Theorem  holds in $\mathfrak X$  for automorphisms if $\AT(G,\phi,H)$ holds for every $G\in \mathfrak X, \ \phi\in \Aut(G)$  and every $\phi$-stable normal subgroup $H$ of $G.$
		\een
	\end{definition}
	We begin with a basic invariance property of the algebraic entropy (see \cite[Lemma 5.1.7]{DG-islam}).
	
	\begin{lemma}\label{lem:iuc}
		Let $G$ and $H$ be groups, $\phi\in \End(G)$ and $\psi\in \End(H)$. If there exists an isomorphism $\xi : G \to H$ such that $\psi=\xi\phi\xi^{-1}$, then
		\[
		h(\phi) = h(\psi).
		\]
	\end{lemma}
	
	We shall also use the following reduction for solvable 
	groups. 
	
	The following reduction was proved in \cite[Proposition~3.3(1)]{S}.
	We include a self-contained proof for convenience.
	
	\begin{prop}[Reduction to the derived subgroup]\label{prop:derived-reduction}
		Let $\mathfrak X$ be a class of solvable groups closed under taking
		subgroups and quotients. Assume that $\AT(G,\phi,G')$ holds for every
		$G\in\mathfrak X$ and every $\phi\in\End(G)$. Then the Addition Theorem
		holds in $\mathfrak X$ for endomorphisms.
	\end{prop}
	
	\begin{proof}
		We prove the assertion by induction on the derived length. More precisely,
		for each $d\ge 1$ we prove that, for every $G\in\mathfrak X$ of derived
		length $d$, every $\phi\in\End(G)$, and every $\phi$-invariant normal
		subgroup $H\unlhd G$, the equality $\AT(G,\phi,H)$ holds. The case $d=1$
		is the abelian Addition Theorem \cite[Theorem~1.1]{DGBabelian}.
		Now let $d>1$ and assume that the assertion holds for all groups in
		$\mathfrak X$ of derived length less than $d$. Let $G\in\mathfrak X$ have
		derived length $d$, let $\phi\in\End(G)$, and let $H\unlhd G$ be
		$\phi$-invariant. Put
		$K=H\cap G'$. Since $\mathfrak X$ is closed under taking subgroups and
		quotients, all groups that occur below belong to $\mathfrak X$. Moreover,
		$G'$, $H$, $H'$, and $K$ are $\phi$-invariant. For brevity, write $h(L)$
		and $h(L/M)$ for the entropies of the corresponding restriction and induced
		endomorphism.
		
		Since $G'$ has derived length strictly smaller than that of $G$, the induction
		hypothesis applied to the normal subgroup $K$ of $G'$ gives
		\begin{equation}\label{eq:derived-Gprime}
		h(G')=h(K)+h(G'/K).
		\end{equation}
		We next prove the analogous equality for $H$. By the assumption of the
		proposition, applied to $H$ and its derived subgroup $H'$, we have
		\begin{equation}\label{eq:derived-Hprime}
		h(H)=h(H')+h(H/H').
		\end{equation}
		Since $H'\le K\le G'$, the group $K$ has derived length strictly smaller than
		that of $G$. Hence the induction hypothesis applied to $K$, with the
		$\phi$-invariant normal subgroup $H'$, yields
		\begin{equation}\label{eq:derived-K}
		h(K)=h(H')+h(K/H').
		\end{equation}
		Finally, $H/H'$ is abelian and $K/H'$ is an invariant subgroup of $H/H'$.
		Therefore the Addition Theorem for abelian groups gives
		\begin{equation}\label{eq:derived-abelian-H}
		h(H/H')=h(K/H')+h(H/K).
		\end{equation}
		Combining \eqref{eq:derived-Hprime}--\eqref{eq:derived-abelian-H}, we obtain
		\begin{equation}\label{eq:derived-H}
		h(H)=h(K)+h(H/K).
		\end{equation}
		
		The assumed equality for $G$ gives
		\begin{equation}\label{eq:derived-G}
		h(G)=h(G')+h(G/G').
		\end{equation}
		Applying the abelian Addition Theorem in $G/G'$ to the invariant subgroup
		$HG'/G'\cong H/K$, we obtain
		\begin{equation}\label{eq:derived-abelian-G}
		h(G/G')=h(H/K)+h(G/HG').
		\end{equation}
		Finally, $G/H\in\mathfrak X$, so the assumption applied to $G/H$ and its
		derived subgroup gives
		\[
		h(G/H)=h((G/H)')+h((G/H)/(G/H)').
		\]
		Using the natural isomorphisms
		\[
		(G/H)'\cong G'/K,
		\qquad
		(G/H)/(G/H)'\cong G/HG',
		\]
		and the invariance of algebraic entropy under conjugation by isomorphisms
		(Lemma~\ref{lem:iuc}), this becomes
		\begin{equation}\label{eq:derived-quotient}
		h(G/H)=h(G'/K)+h(G/HG').
		\end{equation}
		Combining \eqref{eq:derived-Gprime}, \eqref{eq:derived-H},
		\eqref{eq:derived-G}, \eqref{eq:derived-abelian-G}, and
		\eqref{eq:derived-quotient}, we obtain
		\[
		\begin{aligned}
		h(G)
		&=h(K)+h(G'/K)+h(H/K)+h(G/HG')\\
		&=h(H)+h(G/H),
		\end{aligned}
		\]
		so $\AT(G,\phi,H)$ holds.
	\end{proof}
	
	We next record a useful consequence of the Addition Theorem. 
	
	\begin{prop}\label{prop:integer1}
		Let $\mathfrak X$ be a class of torsion solvable groups closed under taking subgroups and quotients. If the Addition Theorem holds in $\mathfrak X$ for endomorphisms, then for every $G\in  \mathfrak X$ and every $\phi\in \End(G)$ either $h(\phi)=\infty$ or $h(\phi)=\log(\alpha)$ for some $\alpha\in \Np.$
	\end{prop}
	
	\begin{proof}
		In case $G$ is abelian (that is, of derived length $1$), the assertion follows from \cite[Proposition 1.3]{DGSZ}. Now suppose that the assertion holds for groups from $\mathfrak X$ of derived length less than $n$, and let $G\in \mathfrak X$ be of derived length $n$ and $\phi\in \End(G)$ with $h(\phi)<\infty$. By our assumption, the Addition Theorem holds in $\mathfrak X$, so
		\[
		h(\phi) = h(\phi \upharpoonright_{G'})+h(\bar \phi),
		\]
		where $\bar\phi\in \End(G/G')$ is the induced endomorphism. Since $G'$ has derived length less than $n$ and $G/G'$ is abelian, the induction hypothesis and the abelian case imply that there exist $\alpha, \beta\in \Np$ such that $h(\phi \upharpoonright_{G'})=\log(\alpha)$ and $h(\bar \phi)=\log(\beta)$. Therefore
		\[
		h(\phi)=\log(\alpha)+\log(\beta)=\log(\alpha\beta),
		\]
		and we are done.
	\end{proof}

	\section{Proving the main result}\label{sec:main}
	
	We first prove the Addition Theorem for invariant central subgroups.
	The next lemma strengthens \cite[Proposition~4.1]{S} by removing the
	torsion assumption.
	
	\begin{lemma}[Superadditivity for central extensions]\label{lem:central-superadditivity}
		Let $G$ be a group, $\phi\in\End(G)$, and $N\le Z(G)$ a $\phi$-invariant subgroup.
		Denote by $\bar\phi$ the endomorphism induced by $\phi$ on $G/N$. Then
		\[
		h(\phi)\ge h(\phi\!\upharpoonright_N)+h(\bar\phi_{G/N}).
		\]
	\end{lemma}
	
	\begin{proof}
		Let $E\in\Pfin{N}$ and $C\in\Pfin{G/N}$, let $\pi:G\to G/N$ be the
		quotient map, and choose $B\in\Pfin{G}$ with $\pi(B)=C$. Put
		\[
		A_n=T_n(\phi\!\upharpoonright_N,E),
		\qquad
		D_n=T_n(\phi,B).
		\]
		Centrality and $\phi$-invariance of $N$ give
		\[
		T_n(\phi,EB)=A_nD_n.
		\]
		For each $q\in\pi(D_n)$ choose $d_q\in D_n$ with $\pi(d_q)=q$.
		For fixed $n$, the sets $A_nd_q$, $q\in\pi(D_n)$, are pairwise disjoint,
		since $A_n\subseteq\ker\pi$. Each has cardinality $|A_n|$, and all are
		contained in $A_nD_n=T_n(\phi,EB)$. Since there are exactly
		$|\pi(D_n)|$ such sets, it follows that
		\[
		|T_n(\phi,EB)|
		\ge |A_n|\,|\pi(D_n)|
		=|T_n(\phi\!\upharpoonright_N,E)|\,|T_n(\bar\phi,C)|.
		\]
		After taking logarithms, dividing by $n$, and passing to the limit, we get
		\[
		H(\phi,EB)
		\ge H(\phi\!\upharpoonright_N,E)+H(\bar\phi,C).
		\]
		Since $E$ and $C$ are arbitrary and $H(\phi,EB)\le h(\phi)$, taking the
		two independent suprema proves the result.
	\end{proof}
	
	\begin{lemma}[Subadditivity for central extensions]\label{lem:central-subadditivity}
		Let $G$ be a group, $\phi\in\End(G)$, and $N\le Z(G)$ a $\phi$-invariant subgroup.
		Denote by $\bar\phi$ the endomorphism induced by $\phi$ on $G/N$. Then
		\[
		h(\phi)\le h(\phi\!\upharpoonright_{N})+h(\bar\phi_{G/N}).
		\]
	\end{lemma}
	
	\begin{proof}
		Fix $F\in\Pfin{G}$, let $\pi:G\to G/N$ be the quotient map, and set
		$Q=\pi(F)$. For $m\in\Np$, put
		\[
		X=T_m(\phi,F),\qquad Y=T_m(\bar\phi,Q),\qquad \psi=\phi^m.
		\]
		Since $\pi(X)=Y$, for each $y\in Y$ choose $s(y)\in X$ such that
		$\pi(s(y))=y$, thereby defining a map $s:Y\to X$. Define
		\[
		A=\{x\bigl(s(\pi(x))\bigr)^{-1}:x\in X\}\subseteq N.
		\]
		Thus every $x\in X$ can be written as $x=a s(y)$ with $a\in A$ and
		$y\in Y$.
		
		\medskip
		\noindent
		\emph{Claim.} For every $k\in\Np$,
		\begin{equation}\label{eq:block-central}
		|T_k(\psi,X)|
		\le |T_k(\psi\!\upharpoonright_N,A)|\,|Y|^k.
		\end{equation}
		Indeed, let
		\[
		S_k=\{s(y_0)\psi(s(y_1))\cdots\psi^{k-1}(s(y_{k-1})):
		y_0,\ldots,y_{k-1}\in Y\}.
		\]
		Then $|S_k|\le |Y|^k$. Write $x_j=a_js(y_j)$, where
		$a_j\in A$ and $y_j\in Y$. Since $\psi^j(a_j)\in N\le Z(G)$, we may
		rearrange the central factors to obtain
		\[
		\begin{aligned}
		x_0\psi(x_1)\cdots\psi^{k-1}(x_{k-1})
		&=a_0\psi(a_1)\cdots\psi^{k-1}(a_{k-1})\\
		&\quad\cdot s(y_0)\psi(s(y_1))\cdots\psi^{k-1}(s(y_{k-1})).
		\end{aligned}
		\]
		Hence
		\[
		T_k(\psi,X)\subseteq T_k(\psi\!\upharpoonright_N,A)S_k,
		\]
		and therefore
		\[
		|T_k(\psi,X)|
		\le |T_k(\psi\!\upharpoonright_N,A)|\,|S_k|
		\le |T_k(\psi\!\upharpoonright_N,A)|\,|Y|^k,
		\]
		which proves the claim.
		
		For $j=0,\ldots,k-1$, since $\psi=\phi^m$ and $X=T_m(\phi,F)$,
		\[
		\psi^j(X)
		=\phi^{mj}(F)\phi^{mj+1}(F)\cdots\phi^{mj+m-1}(F).
		\]
		Multiplying these consecutive blocks gives
		\[
		T_k(\psi,X)=T_{mk}(\phi,F).
		\]
		Consequently,
		\[
		\begin{aligned}
		H(\psi,X)
		&=\lim_{k\to\infty}\frac1k\log|T_{mk}(\phi,F)|\\
		&=m\lim_{k\to\infty}\frac1{mk}\log|T_{mk}(\phi,F)|
		=mH(\phi,F).
		\end{aligned}
		\]
		Taking logarithms in \eqref{eq:block-central}, dividing by $k$, and
		letting $k\to\infty$, we obtain
		\begin{equation}\label{eq:block-entropy}
		mH(\phi,F)
		\le \log|Y|+H(\psi\!\upharpoonright_N,A).
		\end{equation}
		
		Set $B=A\cup\{1\}$. Since $A\subseteq B$ and
		$\psi\!\upharpoonright_N=(\phi\!\upharpoonright_N)^m$,
		\[
		T_k(\psi\!\upharpoonright_N,B)
		\subseteq T_{mk}(\phi\!\upharpoonright_N,B)
		\]
		(the missing factors are filled with the identity). Therefore
		\[
		H(\psi\!\upharpoonright_N,A)
		\le H(\psi\!\upharpoonright_N,B)
		\le mH(\phi\!\upharpoonright_N,B)
		\le mh(\phi\!\upharpoonright_N).
		\]
		Substituting this into \eqref{eq:block-entropy}, dividing by $m$, and
		recalling that $Y=T_m(\bar\phi,Q)$, we get
		\[
		H(\phi,F)
		\le h(\phi\!\upharpoonright_N)
		+\frac1m\log|T_m(\bar\phi,Q)|.
		\]
		Letting $m\to\infty$ and using the definition of $H(\bar\phi,Q)$ yields
		\[
		H(\phi,F)
		\le h(\phi\!\upharpoonright_N)+H(\bar\phi,Q)
		\le h(\phi\!\upharpoonright_N)+h(\bar\phi_{G/N}),
		\]
		Taking the supremum over $F\in\Pfin{G}$ completes the proof.
	\end{proof}

	Combining Lemmas~\ref{lem:central-superadditivity} and
	\ref{lem:central-subadditivity}, we immediately obtain the following.
	
	\begin{corol}[Addition Theorem for central extensions]\label{cor:atforcen}
		Let $G$ be a group, $\phi\in\End(G)$, and $N\le Z(G)$ a
		$\phi$-invariant subgroup. Then
		\[
		h(\phi)= h(\phi\upharpoonright_{N})+h(\bar\phi_{G/N}).
		\]
	\end{corol}
	\begin{prop}[Addition Theorem with respect to the center for automorphisms]\label{prop:center-AT}
	Let $G$ be a group and let $\varphi\in\Aut(G)$. Then
	\[
	h(\varphi)
	=
	h(\varphi\!\upharpoonright_{Z(G)})
	+
	h(\bar\varphi_{G/Z(G)}).
	\]
	In other words, $\AT(G,\varphi,Z(G))$ holds for every automorphism
	of every group.
\end{prop}

\begin{proof}
Since $\varphi$ is an automorphism, $\varphi(Z(G))=Z(G)$.
The claim follows from Corollary~\ref{cor:atforcen}.
\end{proof}
	
	Let $G$ be a nilpotent group of class $n$ and consider its fully characteristic central subgroup $\gamma_n(G).$ By Corollary \ref{cor:atforcen} we obtain the following.
	\begin{corol}\label{cor:gamman}
		Let $G$ be a nilpotent group of class $n$ and $\phi\in \End(G).$  Then 
		\[
		h(\phi)= h(\phi\upharpoonright_{\gamma_n(G)})+h(\bar\phi_{G/\gamma_n(G)}).
		\]
	\end{corol}
	Now we are ready to prove our main result.
	
	\begin{thm}\label{thm:main}
		Let $G$ be a nilpotent group, let $\phi\in\End(G)$ and let $H$ be a $\phi$-invariant normal subgroup of $G$. Then $\AT(G,\phi,H)$ holds.
	\end{thm} 
	\begin{proof}
		For each $n\ge 1$, let $\mathfrak X_n$ be the class of nilpotent groups of
		nilpotency class at most $n$. We prove by induction on $n$ that the Addition
		Theorem holds in $\mathfrak X_n$ for endomorphisms.
		
		For $n=1$, all groups in $\mathfrak X_1$ are abelian, so the result follows
		from the Addition Theorem for abelian groups \cite[Theorem~1.1]{DGBabelian}.
		Assume now that $n>1$ and that the Addition Theorem holds for every group in
		$\mathfrak X_{n-1}$. Since $\mathfrak X_n$ is closed under taking subgroups
		and quotients, Proposition \ref{prop:derived-reduction} shows that it is
		enough to prove $\AT(G,\phi,G')$ for every $G\in\mathfrak X_n$ and every
		$\phi\in\End(G)$.
		
		If $G\in\mathfrak X_{n-1}$, this already follows from the induction
		hypothesis. Hence let $G$ have nilpotency class exactly $n$. Applying
		Corollary \ref{cor:gamman}, we obtain
		\begin{equation}\label{eq:cenandcom0}
			h(\phi)= h(\phi\upharpoonright_{\gamma_n(G)})
			+h(\bar\phi_{G/\gamma_n(G)}).
		\end{equation}
		
		\medskip
		
		\noindent
		\emph{Claim 1.}
		\[
		h(\phi\upharpoonright_{G'})
		=h(\phi\upharpoonright_{\gamma_n(G)})
		+h(\bar\phi_{G'/\gamma_n(G)}).
		\]
		
		\begin{proof}[Proof of Claim 1]
			Since $\gamma_k(G')\le \gamma_{k+1}(G)$ for every $k\ge1$, the subgroup
			$G'$ has nilpotency class at most $n-1$. Moreover,
			$\gamma_n(G)\le G'$ is $\phi$-invariant. Thus the induction hypothesis,
			applied to $G'$ and the invariant normal subgroup $\gamma_n(G)$, gives
			the claimed equality.
		\end{proof}
		
		\medskip
		
		\noindent
		\emph{Claim 2.}
		\[
		h(\bar\phi_{G/\gamma_n(G)})
		= h(\bar\phi_{G'/\gamma_n(G)})+h(\bar\phi_{G/G'}).
		\]
		
		\begin{proof}[Proof of Claim 2]
			The quotient $G/\gamma_n(G)$ is nilpotent of class at most $n-1$, and
			$G'/\gamma_n(G)$ is a $\bar\phi$-invariant normal subgroup of it.
			Therefore the induction hypothesis gives the Addition Theorem for
			$G/\gamma_n(G)$ with respect to $G'/\gamma_n(G)$. Since
			\[
			(G/\gamma_n(G))/(G'/\gamma_n(G))\cong G/G',
			\]
			Lemma \ref{lem:iuc} yields the claimed equality.
		\end{proof}
		
		Combining Claims~1 and~2 with \eqref{eq:cenandcom0}, we obtain
		\begin{align*}
		h(\phi)
		&=h(\phi\upharpoonright_{\gamma_n(G)})
		+h(\bar\phi_{G/\gamma_n(G)})\\
		&=h(\phi\upharpoonright_{\gamma_n(G)})
		+h(\bar\phi_{G'/\gamma_n(G)})+h(\bar\phi_{G/G'})\\
		&=h(\phi\upharpoonright_{G'})+h(\bar\phi_{G/G'}).
		\end{align*}
		Thus $\AT(G,\phi,G')$ holds. Hence it holds for every
		$G\in\mathfrak X_n$, and Proposition \ref{prop:derived-reduction} implies
		that the Addition Theorem holds in $\mathfrak X_n$ for endomorphisms.
		This completes the induction and the proof.
	\end{proof}
	
	It was proved by Giordano Bruno and Spiga \cite[Theorem~3.2(b)]{GBSp} that inner automorphisms
	of arbitrary groups have the same algebraic entropy as the identity.
	Since every finite subset of a nilpotent group generates a finitely generated
	nilpotent subgroup, and such groups have polynomial growth by a classical theorem
	of Wolf \cite[Theorem~3.2]{Wolf}, for every finite $F\subseteq G$ there exist $C,d>0$ such that
	$|F^n|\le Cn^d$. Hence
	\[
	H(\operatorname{id}_G,F)=\lim_{n\to\infty}\frac{1}{n}\log|F^n|=0,
	\]
	and therefore $h(\operatorname{id}_G)=0$. It follows that every inner automorphism
	of a nilpotent group has zero algebraic entropy.
	
	This conclusion can also be recovered directly from Theorem~\ref{thm:main},
	by applying the Addition Theorem to the central extension
	$Z(G)\trianglelefteq G$ and arguing by induction on the nilpotency class.
	\begin{corol}\label{cor:inner-zero}
		Let $G$ be a nilpotent group and let $g\in G$.
		Then the inner automorphism $\iota_g\colon G\to G$,
		$\iota_g(x)=g^{-1}xg$, has zero algebraic entropy.
	\end{corol}
	Using Theorem \ref{thm:main} and Proposition \ref{prop:integer1} we deduce the following result.
	\begin{corol} \label{cor:metanip} If $G$ is a torsion nilpotent  group, and $\phi\in \End(G)$ is of finite entropy,  then  $h(\phi)=\log(\alpha)$ for some $\alpha\in \Np.$
	\end{corol}
	As an application of the Addition Theorem for nilpotent groups,
	we now extend Proposition~\ref{prop:center-AT} from the center to all
	terms of the upper central series.
	\begin{prop}[Addition Theorem for the $n$-center of automorphisms]\label{prop:Zn-center-AT}
		Let $G$ be a group and let $\varphi\in\Aut(G)$.
		Then for every $n\in\Np$,
		\[
		h(\varphi)
		=
		h(\varphi\!\upharpoonright_{Z_n(G)})
		+
		h(\bar\varphi_{G/Z_n(G)}).
		\]
		In other words, $\AT(G,\varphi,Z_n(G))$ holds for every $n\in\Np$.
	\end{prop}
	
	\begin{proof}
		We argue by induction on $n$.
		For $n=1$ the claim is exactly Proposition~\ref{prop:center-AT}.
		
		Assume that the claim holds for some $n\ge 1$, and set $\pi\colon G\to G/Z_n(G)$.
		Let $\bar\varphi$ be the automorphism induced by $\varphi$ on $G/Z_n(G)$.
		Proposition~\ref{prop:center-AT} applied to $(G/Z_n(G),\bar\varphi)$ yields
		\[
		h(\bar\varphi)
		=
		h\!\left(\bar\varphi\!\upharpoonright_{Z(G/Z_n(G))}\right)
		+
		h\!\left(\overline{\bar\varphi}_{(G/Z_n(G))/Z(G/Z_n(G))}\right).
		\]
		Using $Z(G/Z_n(G))=Z_{n+1}(G)/Z_n(G)$ and
		$(G/Z_n(G))/Z(G/Z_n(G))\cong G/Z_{n+1}(G)$, we obtain
		\begin{equation}\label{eq:Zn-quotient}
		h(\bar\varphi_{G/Z_n(G)})
		=
		h(\widetilde\varphi_{Z_{n+1}(G)/Z_n(G)})
		+
		h(\bar\varphi_{G/Z_{n+1}(G)}),
		\end{equation}
		where $\widetilde\varphi$ denotes the automorphism induced by $\varphi$ on
		$Z_{n+1}(G)/Z_n(G)$.
		
		On the other hand, $Z_{n+1}(G)$ is a nilpotent group (of class at most $n+1$),
		and $Z_n(G)\unlhd Z_{n+1}(G)$ is $\varphi$-stable. Therefore,
		Theorem~\ref{thm:main} applied to $\varphi\!\upharpoonright_{Z_{n+1}(G)}$ gives
		\begin{equation}\label{eq:Zn-extension}
		h(\varphi\!\upharpoonright_{Z_{n+1}(G)})
		=
		h(\varphi\!\upharpoonright_{Z_n(G)})
		+
		h(\widetilde\varphi_{Z_{n+1}(G)/Z_n(G)}).
		\end{equation}
		
		Finally, by the inductive hypothesis,
		\begin{equation}\label{eq:Zn-induction}
		h(\varphi)=h(\varphi\!\upharpoonright_{Z_n(G)})+h(\bar\varphi_{G/Z_n(G)}).
		\end{equation}
		Combining \eqref{eq:Zn-quotient}, \eqref{eq:Zn-extension}, and
		\eqref{eq:Zn-induction}, we obtain
		\[
		h(\varphi)
		=
		h(\varphi\!\upharpoonright_{Z_{n+1}(G)})
		+
		h(\bar\varphi_{G/Z_{n+1}(G)}),
		\]
		which completes the induction.
	\end{proof}
	
	\section{Continuity along directed unions and applications}\label{sec:directed-unions}
	We shall need the following basic monotonicity result along directed unions
	(see \cite[Proposition~5.1.10]{DG-islam}).
	\begin{prop}\label{prop:entropy-directed-union}
		Let $G$ be a group and let $\varphi\in\End(G)$.
		Let $\{G_i\}_{i\in I}$ be a \emph{directed} family of $\varphi$-invariant subgroups of $G$
		such that
		\[
		G=\bigcup_{i\in I}G_i.
		\]
		Then
		\[
		h(\varphi)=\sup_{i\in I} h\bigl(\varphi\!\upharpoonright_{G_i}\bigr).
		\]
	\end{prop}
	
	\begin{proof}
		Monotonicity gives
		$\sup_i h(\varphi\!\upharpoonright_{G_i})\le h(\varphi)$.
		Conversely, every $F\in\Pfin{G}$ is contained in some $G_i$, by
		directedness, and then
		\[
		H(\varphi,F)=H(\varphi\!\upharpoonright_{G_i},F)
		\le h(\varphi\!\upharpoonright_{G_i}).
		\]
		Taking the supremum over $F$ gives the reverse inequality.
	\end{proof}
	
	\begin{prop}\label{prop:AT-chain}
		Let $G$ be a group and let $\{G_i\}_{i\in I}$ be a chain of subgroups of $G$
		such that $G=\bigcup_{i\in I}G_i$.
		\begin{enumerate}
			\item\label{prop:AT-chain-endo}
			\emph{(Endomorphisms).}
			Assume that each $G_i$ is \emph{fully characteristic} in $G$ and that $\AT(G_i)$ holds
			for endomorphisms for every $i\in I$. Then $\AT(G)$ holds for endomorphisms.
			
			\item\label{prop:AT-chain-auto}
			\emph{(Automorphisms).}
			Assume that each $G_i$ is \emph{characteristic} in $G$ and that $\AT(G_i)$ holds
			for automorphisms for every $i\in I$. Then $\AT(G)$ holds for automorphisms.
		\end{enumerate}
	\end{prop}
	
	\begin{proof}
		We prove~(\ref{prop:AT-chain-endo}). Let $\varphi\in\End(G)$ and let
		$H\unlhd G$ be $\varphi$-invariant. Since each $G_i$ is fully
		characteristic, $H\cap G_i$ is a $\varphi$-invariant normal subgroup of
		$G_i$, and the Addition Theorem in $G_i$ gives
		\[
		h(\varphi\!\upharpoonright_{G_i})
		=
		h(\varphi\!\upharpoonright_{H\cap G_i})
		+
		h(\bar\varphi\!\upharpoonright_{G_i/(H\cap G_i)}).
		\]
		Identify $G_i/(H\cap G_i)$ with $G_iH/H\le G/H$. The two summands on
		the right form increasing families; hence, since the $G_i$ form a chain,
		the supremum of their sums is the sum of their suprema. Applying
		Proposition~\ref{prop:entropy-directed-union} to
		\[
		G=\bigcup_iG_i,
		\qquad
		H=\bigcup_i(H\cap G_i),
		\qquad
		G/H=\bigcup_iG_iH/H,
		\]
		therefore yields
		\[
		h(\varphi)=h(\varphi\!\upharpoonright_H)+h(\bar\varphi).
		\]
		The automorphism case is identical, using characteristic subgroups in
		place of fully characteristic ones.
	\end{proof}
\begin{corol}\label{cor:AT-Gn!-implies-ATG}
	Let $G$ be a locally finite group.
	If $\AT(G[n!])$ holds for \emph{endomorphisms} for every $n\in\Np$,
	then $\AT(G)$ holds for endomorphisms.
\end{corol}

\begin{proof}
	For a locally finite group $G$, the family $\{G[n!]\}_{n\in\Np}$ forms a chain
	of fully characteristic subgroups whose union is $G$, so the conclusion follows from
	Proposition~\ref{prop:AT-chain}(\ref{prop:AT-chain-endo}).
\end{proof}

\begin{corol}\label{cor:bounded-reduction}
	Let $\mathfrak X$ be a class of locally finite groups closed under taking
	subgroups and quotients. Then the Addition Theorem holds in $\mathfrak X$
	for endomorphisms if and only if it holds for endomorphisms of every group
	in $\mathfrak X$ that is generated by a bounded subset, that is, by a subset
	whose element orders are uniformly bounded.
\end{corol}

\begin{proof}
	One implication is immediate. Conversely, assume that the Addition Theorem
	holds for endomorphisms of every group in $\mathfrak X$ generated by a bounded
	subset, and let $G\in\mathfrak X$. For every $n\in\Np$, the subgroup $G[n!]$
	belongs to $\mathfrak X$ and is generated by the bounded subset
	$\{x\in G:x^{n!}=1\}$. Hence $\AT(G[n!])$ holds for endomorphisms for every
	$n\in\Np$, and Corollary~\ref{cor:AT-Gn!-implies-ATG} yields that $\AT(G)$
	holds for endomorphisms. Since $G$ was arbitrary, the Addition Theorem holds
	in $\mathfrak X$ for endomorphisms.
\end{proof}

\begin{prop}\label{prop:AT-omega-hypercentral}
	Let $G$ be an $\omega$-hypercentral group and let $\phi\in\Aut(G)$.
	Then the Addition Theorem holds for $(G,\phi)$; that is, for every $\phi$-stable
	normal subgroup $H\trianglelefteq G$ one has
	\[
	h(\phi)=h(\phi\!\upharpoonright_H)+h(\bar\phi_{G/H}).
	\]
	Moreover, if $G$ is also torsion, then either $h(\phi)=\infty$ or
	$h(\phi)=\log(\alpha)$ for some $\alpha\in\Np.$
\end{prop}

\begin{proof}
	Since $G$ is $\omega$-hypercentral, we have $G=\bigcup_{n\in\N}Z_n(G)$.
	Each $Z_n(G)$ is a characteristic (hence $\phi$-stable) nilpotent subgroup of
	$G$, of class at most $n$. By Theorem~\ref{thm:main}, $\AT(Z_n(G))$
	holds for endomorphisms, and hence also for automorphisms. Therefore
	Proposition~\ref{prop:AT-chain}(\ref{prop:AT-chain-auto}) yields $\AT(G)$
	for automorphisms.
\medskip
Assume now that $G$ is also torsion. Since every $\omega$-hypercentral group
is locally nilpotent, $G$ is locally finite. Applying
Proposition~\ref{prop:entropy-directed-union} to the directed chain
$\{Z_n(G)\}_{n\in\N}$, we obtain
\[
h(\phi)
=
\sup_{n\in\N}
h\bigl(\phi\!\upharpoonright_{Z_n(G)}\bigr).
\]
By Corollary~\ref{cor:metanip}, for each $n$ either
\[
h\bigl(\phi\!\upharpoonright_{Z_n(G)}\bigr)=\infty
\]
or
\[
h\bigl(\phi\!\upharpoonright_{Z_n(G)}\bigr)
=
\log(\alpha_n)
\qquad
\text{for some }\alpha_n\in\Np.
\]
If $h(\phi)<\infty$, then all these entropies are finite, and hence
\[
h(\phi)
=
\sup_{n\in\N}\log(\alpha_n)
=
\log(\alpha),
\]
where
\[
\alpha=\sup_{n\in\Np}\alpha_n\in\Np.
\]
\end{proof}

We conclude with an example of a torsion $\omega$-hypercentral group
that is not nilpotent, illustrating that the previous proposition
goes beyond the nilpotent case.

\begin{example}
Let $p$ be a prime and consider the restricted direct sum
\[
G=\bigoplus_{m\geq2}\UT_m(p),
\]
where $\UT_m(p)$ denotes the group of upper unitriangular
$m\times m$ matrices over $\mathbb F_p$.

Each $\UT_m(p)$ is a finite $p$-group, and every element of $G$ has
finite support. Hence $G$ is a torsion group.

We next show that $G$ is $\omega$-hypercentral. For every $r\in\N$,
the upper central series is computed coordinatewise:
\[
Z_r(G)
=
\bigoplus_{m\geq2}Z_r\bigl(\UT_m(p)\bigr).
\]
Indeed, this follows by induction on $r$, using
\[
G/Z_r(G)
\cong
\bigoplus_{m\geq2}
\UT_m(p)/Z_r\bigl(\UT_m(p)\bigr).
\]
Since the center of a restricted direct sum is the restricted direct sum of
its centers, the defining relation for the upper central series then gives
the asserted formula for $Z_{r+1}(G)$.

Let $x\in G$. Since $x$ has finite support, there exists $M\geq2$
such that all its non-trivial coordinates occur in factors
$\UT_m(p)$ with $m\leq M$. Since $\UT_m(p)$ has nilpotency class
$m-1$, we have
\[
Z_{M-1}\bigl(\UT_m(p)\bigr)=\UT_m(p)
\qquad
\text{for every }m\leq M.
\]
Consequently, $x\in Z_{M-1}(G)$. Therefore
\[
G=\bigcup_{r\in\N}Z_r(G),
\]
so $G$ is $\omega$-hypercentral.

On the other hand, $G$ is not nilpotent. Indeed, if $G$ had nilpotency
class $c$, then every subgroup of $G$ would have nilpotency class at most
$c$. This contradicts the fact that $G$ contains subgroups isomorphic to
$\UT_m(p)$ for every $m\geq2$, whose nilpotency classes $m-1$ are unbounded.

In particular, by
Proposition~\ref{prop:AT-omega-hypercentral},
$\AT(G,\varphi,H)$ holds for every $\varphi\in\Aut(G)$ and every
$\varphi$-stable normal subgroup $H$ of $G$. Moreover, for every
$\varphi\in\Aut(G)$, either $h(\varphi)=\infty$ or
$h(\varphi)=\log(\alpha)$ for some $\alpha\in\Np$.
\end{example}

\end{document}